# Log-sum-exp optimization problem subjected to Lukasiewicz fuzzy relational inequalities

Amin Ghodousian[*], Alireza Norouzi Azad, Hadi Amiri

*Faculty of Engineering Science, College of Engineering*
*University of Tehran, P.O.Box 11365-4563, Tehran, Iran*
*E-mail: a.ghodousian@ut.ac.ir*

*Department of Engineering Science, College of Engineering,*
*University of Tehran, P.O.Box 11365-4563, Tehran, Iran*
*E-mail: Alireza.Norouzi@ut.ac.ir*

*Faculty of Engineering Science, College of Engineering,*
*University of Tehran, Tehran, Iran*
*E-mail: Hadi.Amiri@ut.ac.ir*

**Abstract**

*In this paper, we introduce a nonlinear optimization problem whose objective function is the convex log-sum-exp function and the feasible region is defined as a system of fuzzy relational inequalities (FRI) defined by the Lukasiewicz t-norm. Some necessary and sufficient conditions are derived to determine the feasibility of the problem. The feasible solution set is characterized in terms of a finite number of closed convex cells. Since the feasible solutions set of FRIs is non-convex, conventional methods may not be directly employed. An algorithm is presented for solving this nonlinear problem. It is proved that the algorithm can find the exact optimal solution and an example is presented to illustrate the proposed algorithm.*

## Introduction

Resolution of fuzzy relational equations (FRE) with max-min composition was first studied by Sanchez [25]. Besides, Sanchez developed the application of FRE in medical diagnosis in biotechnology. Nowadays, it is well known that many issues associated with a body knowledge can be treated as FRE problems [22]. Since then the composition operator in FREs was replaced by max-product and furthermore extended to the general max-t-norm composition operator [3,4,7,8,10,11,17,18,20]. The resolution method was kept in improving. In fact, if the max-t-norm fuzzy relational equations is consistent, then its solution set is often a non-convex set that is completely determined by a unique maximum solution and a finite number of minimal solutions [5,10,11,21,23,24].

Most of the existing literatures in this research field focused on the resolution of FRE and its relevant optimization problems. However, only a few research works investigated the fuzzy relational inequalities (FRI) and its relevant optimization problems [7-9,12,13,16,28]. For example, Guo et al. [12] studied the linear programming problem with max-min FRI constraint. Li and Yang [16] introduced the so-called addition-min FRI to characterize a peer-to-peer file sharing system. Based on the concept of pseudo-minimal index, Yang [27] developed a pseudo-minimal-index algorithm to minimize a linear objective function with addition-min FRI constraint. To improve the results presented in [27], Yang et al. [29] proposed the min-max programming subject to addition-min fuzzy relational inequalities. They also studied the multi-level linear programming problem with addition-min FRI constraint [30]. Although the FRE and FRI with max-min and max-product operators were widely applied, the max-Lukasiewicz operator was superior to the other operators in many specific cases. In [11], the authors solved a nonlinear optimization problem subjected to FREs defined with Lukasiewicz t-norm. The optimization problems with general nonlinear objective functions and FRE or FRI constraints were studied in [6-11,19]. In general, the genetic algorithm was applied to deal with this kind of problems. However, in our previous works, we concentrated on finding the optimum solutions to the linear optimization problems subjected to Lukasiewicz and Aczel-Alsina Fuzzy relational equations [32,33]. also, some fuzzy relation nonlinear optimization problems could be solved by some specific method. For example, fuzzy relation geometric programming problem was introduced by Yang and Cao [2]. Yang et al. [31] studied the single-variable term semi-latticized geometric programming subject to max-product fuzzy relation equations. The proposed problem was devised from the peer-to-peer network system and the target was to minimize the biggest dissatisfaction degrees of the terminals in such system. Yang et al. [28] introduced another version of the latticized programming problem subject to max-prod fuzzy relation inequalities with application in the optimization management model of wireless communication emission base stations. The latticized problem was defined by minimizing objective function $z(x) = \max\{x_1, x_2, ..., x_n\}$ subject to feasible region $X(A,b) = \{x \in [0,1]^n : A \circ x \geq b\}$ where "$\circ$" denotes fuzzy max-product composition. However, in many applications (especially, in convex programing problems), the function $f(x) = \log(e^{x_1} + e^{x_2} + \cdots + e^{x_n})$ (called log-sum-exp function) is used as a differentiable approximation of the max function $\max\{x_1, x_2, ..., x_n\}$ [1].

In this paper, we study the following non-linear optimization problem in which the objective function is defined as the log-sum-exp function and the constraints are formed as the Lukasiewicz fuzzy relational inequalities:

$$\min \log(e^{x_1} + e^{x_2} + \cdots + e^{x_n})$$
$$A \circ x \geq b \quad (1)$$
$$x \in [0,1]^n$$

where $I = \{1, 2, ..., m\}$, $J = \{1, 2, ..., n\}$, $A = (a_{ij})_{m \times n}$, $0 \leq a_{ij} \leq 1$ ($\forall i \in I$ and $\forall j \in J$), is a fuzzy matrix, $b = (b_i)_{m \times 1}$, $0 \leq b_i \leq 1$ ($\forall i \in I$), is an $m$-dimensional fuzzy vector, and "$\circ$" is the max-Lukasiewicz composition, that is, $x \circ y = \max\{x + y - 1, 0\}$. If $a_i$ is the $i$'th row of matrix $A$, then the constraints can be expressed as $a_i \circ x \geq b_i$ ($i \in I$), where $a_i \circ x = \max_{j \in J}\{a_{ij} \circ x_j\}$ and $a_{ij} \circ x_j = \max\{a_{ij} + x_j - 1, 0\}$.

The rest of the paper is organized as follows. In Section 2, some conditions are derived to determine the feasibility of max-Lukasiewicz FRI. In Section 3, the feasible solution set of problem (1) is characterized. It is shown that the feasible region can be expressed as the union of a finite closed convex cells. Section 4 describes the optimal solution of Problem (1). An algorithm is proposed to find the optimal solution and finally, Section 5 provides a numerical example to illustrate the algorithm.

**Feasibility conditions for problem (1)**

In [8], the authors discussed some properties of FRIs, where the FRIs defined by operators with (closed) convex solutions. in this section we introduce some relevant results about the solution to system of max-Lukasiewicz FRI. denote $S(A,b) = \{x \in [0,1]^n : A \circ x \geq b\}$, that is, set $S(A,b)$ represents the solution set of system (1).

Similar to the general mathematical programming problem, a vector $x$ satisfying the constraints in (1), i.e., $x \in S(A,b)$, is said to be a feasible solution of problem (1). According to the foregoing notation, the following result is trivially attained.

**Proposition 1.** If $x \in S(A,b)$, then $x' \in S(A,b)$ holds for any $x' \in [0,1]^n$ satisfying $x' \geq x$.

**Definition 1.** A solution $\bar{x} \in S(A,b)$ is said to be the maximum solution of $S(A,b)$ when $x \leq \bar{x}$ for all $x \in S(A,b)$. A solution $\underline{x} \in S(A,b)$ is said to be a minimal solution of $S(A,b)$ when $x \leq \bar{x}$ implies $x = \bar{x}$ for any $x \in S(A,b)$.

The consistency of problem (1) can be checked by the following theorem.

**Theorem 1.** For each $i \in I$, let $J(i) = \{j \in J : a_{ij} \geq b_i\}$. Then, the following statements are equivalent:
(a) $S(A,b) \neq \emptyset$.
(b) $\mathbf{1} \in S(A,b)$, where $\mathbf{1}$ is an $n$-dimensional vector with each component equal to one.
(c) $J(i) \neq \emptyset$, $\forall i \in I$.

**Proof.** (a) $\Rightarrow$ (b) Let $x \in S(A,b)$. Since $x \leq \mathbf{1}$, then by proposition 1, $\mathbf{1} \in S(A,b)$.
(b) $\Rightarrow$ (c) Suppose that $\mathbf{1} \in S(A,b)$. Therefore, $a_i \circ \mathbf{1} = \max_{j \in J}\{\max\{a_{ij} + 1 - 1, 0\}\} = \max_{j \in J}\{a_{ij}\} \geq b_i$, $\forall i \in I$. Hence, there exists at least one $j_i \in J$ such that $a_{ij_i} \geq b_i$.
(c) $\Rightarrow$ (a) Let $J(i) \neq \emptyset$, $\forall i \in I$. Thus, for each $i \in I$, there exists at least one $j_i \in J$ such that $a_{ij_i} \geq b_i$. So, similar to part (a), we have $a_i \circ \mathbf{1} = \max_{j \in J}\{a_{ij}\} \geq a_{ij_i} \geq b_i$, $\forall i \in I$. Thus, $\mathbf{1} \in S(A,b)$. $\square$

**Corollary 1.** Suppose that $S(A,b) \neq \emptyset$. Then, vector $\mathbf{1}$ is the maximum solution of $S(A,b)$.

**Characterization of the feasible solution set**

Let $S_i(A,b) = \{x \in [0,1]^n : a_i \circ x \geq b_i\}$, $\forall i \in I$. So, it is clear that $S(A,b) = \bigcap_{i \in I} S_i(A,b)$. Based on Corollary 1, vector $\mathbf{1}$ is the maximum solution of $S_i(A,b)$, $\forall i \in I$. The following theorem shows that the set $S_i(A,b)$ has exactly $|J(i)|$ minimal solutions, where $|J(i)|$ denotes the cardinality of the set $J(i)$.

**Theorem 2.** Suppose that $S(A,b) \neq \emptyset$, $i \in I$ and $j_0 \in J(i)$. Also, define $\underline{x}(i, j_0) \in [0,1]^n$ where

$$\underline{x}(i, j_0)_j = \begin{cases} 1 + b_i - a_{ij_0} &, j = j_0 \\ 0 &, otherwise \end{cases}$$

Then, $\underline{x}(i,j_0)$ is a minimal solution for $S_i(A,b)$.

**Proof.** By contradiction, suppose that $x' \in S_i(A,b)$, $x' \leq \underline{x}(i,j_0)$ and $x' \neq \underline{x}(i,j_0)$. So, $x'_j \leq \underline{x}(i,j_0)_j$, $\forall j \in J$ and $x' \neq \underline{x}(i,j_0)$. Therefore, $x'_j = 0$, $\forall j \in J - \{j_0\}$, and $x'_{j_0} < 1 + b_i - a_{ij_0}$. Hence, we have $a_{ij} \circ x'_j = \max\{a_{ij}-1,0\} = 0$, $\forall j \in J - \{j_0\}$, and $a_{ij_0} \circ x'_{j_0} < b_i$. Therefore, $a_i \circ x' < b_i$ which contradicts $x' \in S_i(A,b)$. □

**Corollary 2.** $S_i(A,b) = \bigcup_{j \in J(i)} [\underline{x}(i,j),1]$, $\forall i \in I$.

**Definition 2.** Let $e: I \to \bigcup_{i \in I} J(i)$ so that $e(i) \in J(i)$, $\forall i \in I$, and let $E$ be the set of all vectors $e$. Also, for each $e \in E$, we define $\underline{x}(e) \in [0,1]^n$ such that $\underline{x}(e)_j = \max_{i \in I}\{\underline{x}(i,e(i))_j\}$, $\forall j \in J$.

Theorem 3 below shows that the feasible solutions set of problem (1) can be written in terms of a finite closed convex cells.

**Theorem 3.** $S(A,b) = \bigcup_{e \in E} [\underline{x}(e),1]$.

**Proof.** From Corollary 2 and the equality $S(A,b) = \bigcap_{i \in I} S_i(A,b)$, we have $S(A,b) = \bigcap_{i \in I} \bigcup_{j \in J(i)} [\underline{x}(i,j),1]$, or equivalently $S(A,b) = \bigcup_{e \in E} \bigcap_{i \in I} [\underline{x}(i,e(i)),1]$. Therefore, $S(A,b) = \bigcup_{e \in E} [\max_{i \in I} \{\underline{x}(i,e(i))\},1]$. Now, the result follows from the definition of $\underline{x}(e)$. □

## Resolution of problem (1)

Based on Theorem 3, for each $x' \in S(A,b)$ there exist $e' \in E$ such that $x' \in [\underline{x}(e'),1]$. In other words, if $x' \leq \underline{x}(e)$ and $x' \neq \underline{x}(e)$, $\forall e \in E$, then $x' \notin S(A,b)$. So, we have the following corollary.

**Corollary 3.** Let $\underline{S}(A,b)$ denote the set of all minimal solutions (minimal solution set) of $S(A,b)$. Then, $\underline{S}(A,b) \subseteq \{\underline{x}(e) : e \in E\}$.

The following theorem characterizes the optimal solution of problem (1).

**Theorem 4.** If $S(A,b) \neq \emptyset$, there exists a minimal solution of $S(A,b)$, i.e. $\underline{x}^* \in \underline{S}(A,b)$, such that $\underline{x}^*$ is an optimal solution of problem (1).

**Proof.** Let $f(x) = \log(\sum_{j=1}^n e^{x_j})$. Furthermore, suppose that $M = \min\{f(\underline{x}) : \underline{x} \in \underline{S}(A,b)\}$ and $f(\underline{x}^*) = M$ where $\underline{x}^*$ is a minimal solution. Based on Theorem 3 and Corollary 3, for each $x' \in S(A,b)$ there exist some $\underline{x} \in \underline{S}(A,b)$ such that $\underline{x} \leq x'$, i.e., $\underline{x}_j \leq x'_j$, $\forall j \in J$. So, since $f(x)$ is an increasing function, we have $f(\underline{x}) \leq f(x')$. But, $f(\underline{x}^*) \leq f(\underline{x})$ which implies $f(\underline{x}^*) \leq f(x')$, $\forall x' \in S(A,b)$. □

By combination of Theorem 4 and Corollary 3, it turns out that the optimal solution of problem (1) must be a vector $\underline{x}(e^*)$ for some $e^* \in E$. Based on this fact, we can find the optimal solution of problem (1) by pairwise comparison between the elements of set $\{\underline{x}(e) : e \in E\}$. We now summarize the preceding discussion as an algorithm.

### Algorithm 1
Given problem (1):
**1.** If $\mathbf{1} \notin S(A,b)$, then stop; $S(A,b)$ is empty (Theorem 1).
**2.** Find solutions $\underline{x}(e), \forall e \in E$ (Definition 2).
**3.** Find the minimal solutions, $\underline{S}(A,b)$ by the pairwise comparison between the solutions $\underline{x}(e)$ (Corollary 3).
**4.** Find the optimal solution $\underline{x}(e^*)$ for problem (1) by the pairwise comparison between the objective values of the elements of $\underline{S}(A,b)$ (Theorem 4).

### Numerical example

**Example 1.** Consider the following linear optimization problem (1):

$$\min \ f(x) = \log(\sum_{j=1}^7 e^{x_j})$$

$$\begin{bmatrix} 0.8147 & 0.0975 & 0.1576 & 0.1418 & 0.6557 & 0.7577 & 0.7060 \\ 0.2784 & 0.9058 & 0.9705 & 0.4217 & 0.0357 & 0.7431 & 0.0318 \\ 0.1270 & 0.5468 & 0.9571 & 0.9157 & 0.8491 & 0.3922 & 0.2769 \\ 0.5134 & 0.3575 & 0.4853 & 0.7922 & 0.9339 & 0.6554 & 0.0461 \\ 0.6323 & 0.4648 & 0.8002 & 0.6594 & 0.6787 & 0.1711 & 0.0971 \end{bmatrix} \circ x \geq \begin{bmatrix} 0.7898 \\ 0.8456 \\ 0.9463 \\ 0.7094 \\ 0.7547 \end{bmatrix}$$

$$x \in [0,1]^7$$

where "$\circ$" is the max-Lukasiewicz composition.

**Step 1:** By a simple calculation, we have

$$A \circ \mathbf{1} = \begin{bmatrix} 0.8147 \\ 0.9705 \\ 0.9571 \\ 0.9339 \\ 0.8002 \end{bmatrix} \geq \begin{bmatrix} 0.7898 \\ 0.8456 \\ 0.9463 \\ 0.7094 \\ 0.7547 \end{bmatrix} = b$$

Therefore, the problem is feasible. On the other hand, in this example, $J(1) = \{1\}$, $J(2) = \{2,3\}$, $J(3) = \{3\}$, $J(4) = \{4,5\}$ and $J_5 = \{3\}$. So, from Theorem 1, $J(i) \neq \emptyset$, $\forall i \in \{1,...,5\}$.

**Step 2:** According to Definition 2, $|E|=\prod_{i=1}^{5}|J(i)|=4$. Therefore, this problem has four minimal candidate solutions $\underline{x}(e)$. These solutions are attained as follows:

$e_1 = [1, 2, 3, 4, 3]$
$\underline{x}(e_1) = [0.9751, 0.9398, 0.9892, 0.9172, 0, 0, 0]$

$e_2 = [1, 3, 3, 4, 3]$
$\underline{x}(e_2) = [0.9751, 0, 0.9892, 0.9172, 0, 0, 0]$

$e_3 = [1, 2, 3, 5, 3]$
$\underline{x}(e_3) = [0.9751, 0.9398, 0.9892, 0, 0.7755, 0, 0]$

$e_4 = [1, 3, 3, 5, 3]$
$\underline{x}(e_4) = [0.9751, 0, 0.9892, 0, 0.7755, 0, 0]$

where $e=[j_1, j_2, ..., j_m]$ means $e(i)=j_i \in J(i)$, $\forall i \in I$.

**Step 3:** It is easy to see that $\underline{x}(e_2) \leq \underline{x}(e_1)$ and $\underline{x}(e_4) \leq \underline{x}(e_3)$. Therefore, $\underline{x}(e_2)$ and $\underline{x}(e_4)$ are the only minimal solutions of the problem. Additionally, from Theorem 3 we have $S(A,b) = [\underline{x}(e_2), \mathbf{1}] \cup [\underline{x}(e_4), \mathbf{1}]$. So, the feasible solutions set of the problem is formed as the union of two closed convex cells.

**Step 4:** By a simple pairwise comparison between the objective values of the $\underline{x}(e_2)$ and $\underline{x}(e_4)$, we have

$2.4434 = f(\underline{x}(e_4)) \leq f(\underline{x}(e_2)) = 2.4717$

Therefore, $e^* = e_4$ and $\underline{x}(e_4)$ is the optimal solution of the problem.

**Conclusion**

Considering the practical applications of the max-Lukasiewicz fuzzy relational inequalities in FRI theory and that of the log-sum-exp function in the convex optimization, a nonlinear optimization problem was studied with log-sum-exp as the objective function subjected to the Lukasiewicz FRI. Since a system of the Lukasiewicz FRI is a non-convex set, an algorithm was presented to find an optimal solution by using the structural properties of the problem. For this purpose, some feasibility conditions were firstly derived and then, the feasible region was completely determined in terms of one maximum and a finite number of minimal solutions. It is proved that we can find the exact optimal solution of the proposed problem from the minimal solutions of the constraints, i.e., a system of max-Lukasiewicz FRI. Additionally, a numerical example was given to illustrate the presented algorithm.